\documentclass[english,12pt]{amsart}

\usepackage{preprint}
\usepackage{graphicx}

\begin{document}


\title{Stagnation zones for $\A$-harmonic functions on canonical
domains}

\author[V.M. Miklyukov, A. Rasila and M. Vuorinen]
{Vladimir M. Miklyukov, Antti Rasila and Matti Vuorinen}
\address[Vladimir M. Miklyukov]
{Department of Mathematics, Volgograd State University,
Universitetskii prospect 100, Volgograd 400062, RUSSIA,
Fax + tel: +7-8442 471608}\email{miklyuk@mail.ru}\
\address[Antti Rasila]
{Department of Mathematics and Systems Analysis, Aalto University, P.O. Box 11100,
FI-00076 Aalto, FINLAND, Fax +358-9-451 3016}
\email{antti.rasila@iki.fi}\
\address[Matti Vuorinen]
{Department of Mathematics, FI-20014 University of Turku,
FINLAND} \email{vuorinen@utu.fi}






 \begin{abstract}
We study stagnation zones of $\A$-harmonic
functions on canonical domains in the Euclidean $n$-dimensional space.
Phragm\'en-Lindel\"of type theorems are proved.
 \end{abstract}


\maketitle



\newcounter{minutes}\setcounter{minutes}{\time}
\divide\time by 60
\newcounter{hours}\setcounter{hours}{\time}
\multiply\time by 60
\addtocounter{minutes}{-\time}

\maketitle

\begin{center}
\texttt{FILE:~\jobname .tex,
         printed: \number\year-\number\month-\number\day,
         \thehours.\ifnum\theminutes<10{0}\fi\theminutes}
\end{center}

\section{Introduction}
In this article we investigate solutions of the $\A$-Laplace equation
on canonical domains in the $n$-dimensional Euclidean space.

Suppose that $D$ is a domain in $\Rn$, and let $f\colon D\to {\R}$ be a 
function. For $s>0$, a subset  $\Delta\subset D$ is called {\it 
$s$-zone} ({\it stagnation zone with the deviation $s$})
of $f$, if there exists a constant $C$ such that the difference between
$C$ and the function $f$ is smaller than $s$ on $\Delta$. We may, for
example, consider difference in the sense of the sup norm,
$$
\|f(x)-C\|_{C(\Delta)} = \sup_{x\in \Delta} |f(x)-C| < s\,,
$$
the $L^p$-norm
$$
\|f(x)-C\|_{L^p(\Delta)}=
\bigg(\int\limits_{\Delta}|f(x)-C|^p\,d{\mathcal
H}^n\bigg)^{1/p} < s\,,
$$
or the Sobolev norm
$$
\|f(x)-C\|_{W^1_p(\Delta)}=
\bigg(\int\limits_{\Delta}|\nabla
f(x)|^pd{\mathcal H}^n\bigg)^{1/p}<s\,,
$$
where ${\mathcal H}^d$ is the $d$-dimensional Hausdorff measure in
$\R^n$.

For discussion about the history of the question, recent results and 
applications, see \cite{SS06,SS07}.

Some estimates of stagnation zone sizes for solutions of
the $\mathcal A$-Laplace equation on locally Lipschitz surfaces and 
behaviour of solutions in stagnation zones, were given in \cite{Mik07}. 
In this paper we consider solutions of the $\A$-Laplace equation in 
subdomains of $\R^n$ of a special form, canonical domains. In two-dimensional case, such 
domains are sectors and strips. In higher dimensions, they are 
conical and cylindrical regions. The special form of domains allows us 
to obtain more precise results.

Below we study stagnation zones
of generalized solutions of the $\A$-Laplace equation
$$
{\rm div}\,\A(x,\nabla f)=0
$$
(see \cite{[5]})  with boundary conditions of types (see
definitions \ref{def1} and \ref{def2} below):
$$
\langle \A(x,\nabla f), \overline{\bf n}\rangle =0\,,\quad x\in
\partial D\setminus G
$$
and
$$
f\,\langle \A(x,\nabla f), \overline{\bf n}\rangle =0\,,
\quad x\in \partial D\setminus G
$$
on canonical domains in the Euclidean $n$-dimensional space, where $G$ 
is a closed subset of $\partial D$. We will prove
Phragm\'en-Lindel\"of type theorems for solutions of the
$\A$-Laplace equation with such boundary conditions.

\subsection*{Canonical domains}{}
Let $n\geq 2$. Fix an integer $k$, $1\le k\le n,$ and set
$$
d_k(x)=\Big(\sum\limits_{i=1}^kx_i^2\Big)^{1/2}\,.
$$
We call the set
$$
B_k(t)=\{x\in {\mathbb{R}}^n:d_k(x)<t\}
$$
a $k$-ball and
$$
\Sigma_k(t)=\{x\in {\mathbb{R}}^n:d_k(x)=t\}
$$
a $k$-sphere in ${\mathbb{R}}^n$.
In particular, the symbol $\Sigma_k(0)$ denotes the $k$-sphere
with the radius $0$, i.e. the set
$$
\Sigma_k(0)=\big\{x=(x_1,\ldots,x_k,\ldots,x_n):x_1=\ldots=x_k=0\big\}.
$$
For every $0< k <n$ we set
$$
p_k(x)=\Big(\sum_{j=k+1}^nx_j^2\Big)^{1/2}\,.
$$
and
$$
\Sigma^*_k(t)=\{x\in {\mathbb{R}}^n:p_k(x)=t\}, \qquad t\geq 0.
$$

Let $0<\alpha<\beta<\infty$ be fixed, and let
$$
D^k_{\alpha,\beta}=\{x\in \Rn: \alpha<p_k(x)<\beta\}.
$$
For $k=n-1$ we also assume that $x_n>0$.
Then for $k=n-1$ the $D^{n-1}_{\alpha,\beta}$ is the a layer between two
parallel hyperplanes, and for $1\le k <n-1$ the boundary of the domain
$D^k_{\alpha,\beta}$  consists of two coaxial cylindrical surfaces. The
intersections $\Sigma_k(t)\cap D^k_{\alpha,\beta}$ are precompact for
all $t>0$. Thus, the functions $d_k(x)$ are {\it exhaustion functions} 
for $D^k_{\alpha,\beta}$.

\begin{center}
\includegraphics[width=6cm]{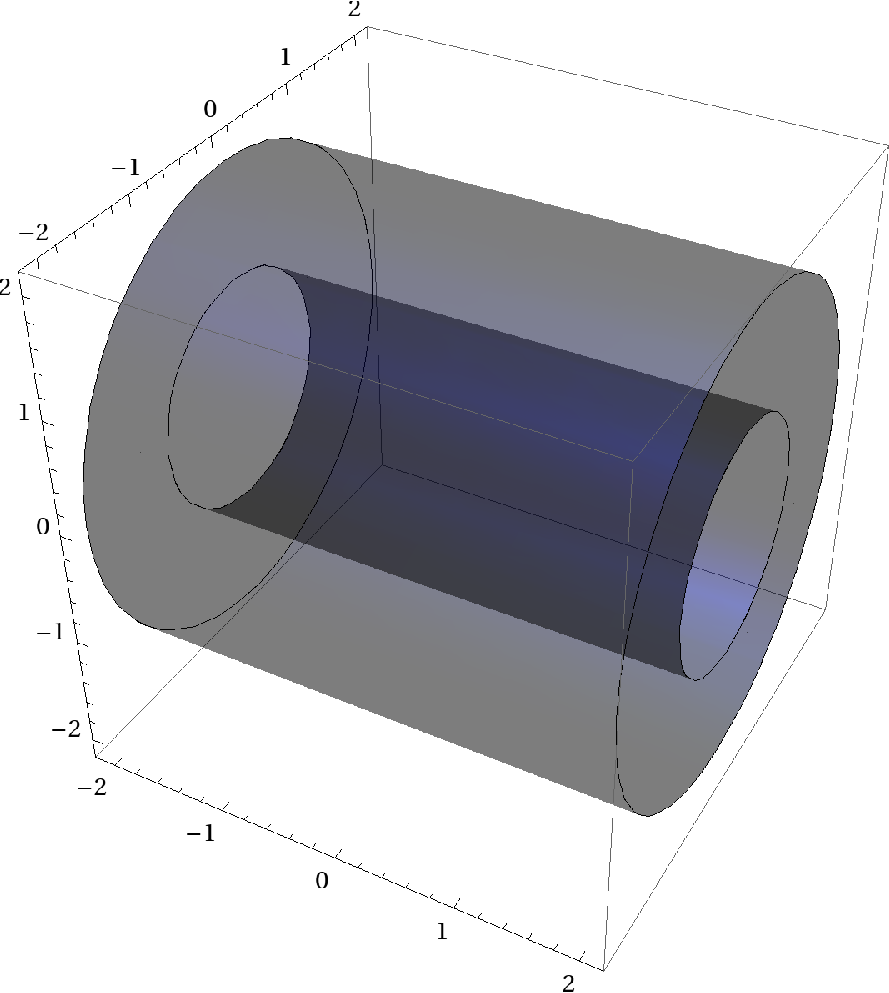}
\qquad 
\includegraphics[width=6cm]{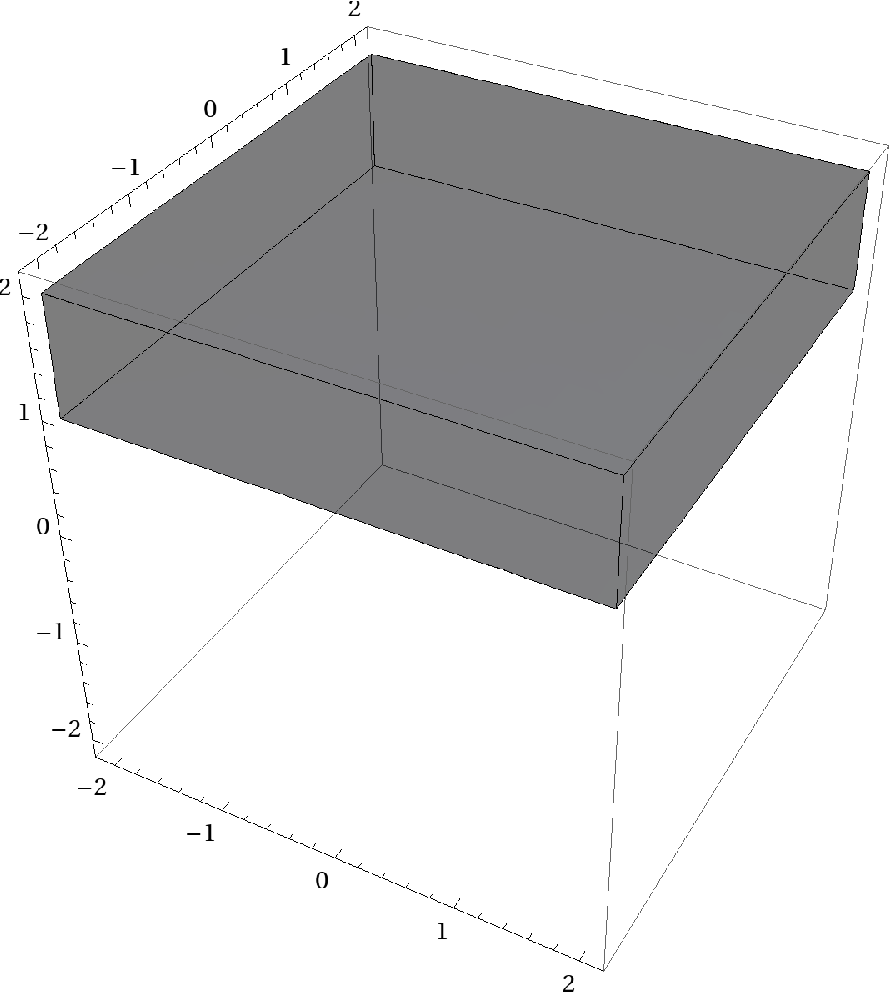}
\end{center}

\noindent
{\sc Figure:} $D^1_{1,2}$ (left) and $D^2_{1,2}$ in $\R^3$.

\medskip

\subsection*{Structure conditions}{}
Let $D$ be a subdomain of $\Rn$ and let
$$
\A(x,\xi)\colon \overline{D}\times {\mathbb{R}}^n \to {\mathbb{R}}^n
$$
be a vector function such that for a.e. $x\in \overline{D}$
the function
$$
\A(x,\xi )\colon {\mathbb{R}}^n \to {\mathbb{R}}^n
$$
is defined and continuous with respect to $\xi$. We assume that the 
function
$$
x \mapsto \A(x,\xi)
$$
is measurable in the Lebesgue sense for all 
$\xi\in \Rn$ and
\begin{equation}
\label{areq5}
\A(x,\lambda \xi)=\lambda\,|\lambda|^{p-2}\,\A(x,\xi)\,,\quad
\lambda\in {\mathbb{R}\setminus\{0\}},\; p\ge 1.
\end{equation}
Suppose that for a.e. $x\in D$ and for all $\xi \in \Rn$
the following properties hold:
\begin{equation}
\label{eq2.23}
\nu_1\,|\xi|^{p}\le \langle\xi, \A(x,\xi)\rangle,
\quad
|\A(x,\xi)|\le \nu_2\,|\xi|^{{p}-1}
\end{equation}
with ${p}\ge 1$ and some constants $\nu_1,\,\nu_2>0$.
We consider the equation
\begin{equation}
\label{eq2.19} {\rm div}\,
\A(x,\nabla f)=0\;.
\end{equation}
An important special case of (\ref{eq2.19}) is the Laplace equation
$$
\Delta f=
\sum_{i=1}^n{\frac{\partial^2 f}{\partial x_i^2}}=0\,.
$$
As in \cite[Chapter 6]{[5]}, we call continuous weak solutions of the
equation (\ref{eq2.19}) {\it $\A$-harmonic} functions. However we should
note that our definition of generalized solutions is slightly
different from the definition given in \cite[p.56]{[5]}.

\subsection*{Frequencies}{}
Fix $t\ge 0$ and $p\ge 1$. Let $O$ be an open subset of $\Sigma^*_k(t)$
(with respect to the relative topology of $\Sigma^*_k(t)$), and let
$\mathcal P$ be a nonempty closed subset of $\partial O$. We set
\begin{equation}
\label{eqla1} \lambda_{p,{\mathcal P}}(O)=\inf\limits_{u}\frac{
\displaystyle\int\limits_O |\nabla u|^p d{\mathcal H}^{n-1}}
{\displaystyle\int\limits_O  u^p d{\mathcal H}^{n-1}}\,,
\end{equation}
where $u\in{\rm Lip}_{\rm loc}(O)\cap C^0(\overline{O})$ with $u|_{{\mathcal P}}=0$.
If ${\mathcal P}=\partial O$ we call  $\lambda_p(O)\equiv\lambda_{p,{\mathcal
P}}(O)$ the
{\it first frequency} of the order $p\ge 1$ of the set $O$.  If ${\mathcal
P}\ne
\partial O$ the quantity $\lambda_{p,{\mathcal P}}(O)$ is the {\it third
frequency}.

The {\it second frequency} is the following quantity:
\begin{equation}
\label{eqla2}
\mu_{p}(O)=\sup\limits_C\inf\limits_{u}{
\frac{\displaystyle\int\limits_O |\nabla u|^p d{\mathcal H}^{n-1}}
{\displaystyle\int\limits_O  (u-C)^p d{\mathcal H}^{n-1}}}\,,
\end{equation}
where the supremum is taken over all constants $C$ and
$u\in {\rm Lip}_{\rm loc}(O)\cap C^0(\overline{O})\,$.
See also P\'olya and Szeg\"o \cite{PS}, Lax \cite{lax}.

\subsection*{Generalized boundary conditions}
Suppose that $D$ is a proper subdomain of ${\mathbb{R}}^n$. Let $\varphi\colon D \to{\mathbb{R}}$ be a locally  Lipschitz function. We denote by $D_b(\varphi)$ the set of all points $x\in D$ at which $\varphi$ does not have the differential.
Let $U\subset D$ be a subset and let $\partial' U=\partial U\setminus \partial D$
be its boundary with respect to $D$. If $\partial' U$ is
$({\mathcal H}^{n-1},n-1)$-rectifiable, then it has locally finite
perimeter in the sense of De~Giorgi, and therefore a unit normal vector $\overline{\mathbf{n}}$ exists $\mathcal{H}^{n-1}$-almost everywhere on $\partial U$ \cite[Sections 3.2.14, 3.2.15]{Fe}.

Let $D\subset\Rn$ be a domain and let $G\subset \partial D$
be a subset of the boundary of $D$. Define the concept of a generalized 
solution of (\ref{eq2.19}) with  zero boundary conditions on $\partial 
D\setminus G$. A subset $U\subset D$ is called {\it admissible}, if 
$\overline{U}\cap\overline{G}=\emptyset$ and $U$ has a $({\mathcal 
H}^{n-1},n-1)$-rectifiable boundary with respect to $D$.

Suppose that $D$ is unbounded.  Let  $G\subset {\partial} D$
be a set closed in $\Rn\cup\{\infty\}$. We denote by $(G,D)$ the collection of 
all subdomains $U\subset D$ with ${\partial} U\subset \left(D\cup 
(\partial D\setminus G)\right)$ and $({\mathcal H}^{n-1},n-1)$-rectifiable
boundaries $\partial' U={\partial} U\setminus {\partial} D$.

\begin{defn}
\label{def1}
We say that a locally Lipschitz function $f\colon D\to{\mathbb{R}}$ is a generalized solution of (\ref{eq2.19}) with the boundary condition
\begin{equation}
\label{areq3}
\langle \A(x,\nabla f), \overline{\bf n}\rangle =0\,,\quad x\in
\partial D\setminus G\,,
\end{equation}
if for every subdomain $U\in (G,D)$,
\begin{equation}
\label{eqphiU}
{\mathcal H}^{n-1}\big[\partial' U \cap D_b(f) \big]=0\,,
\end{equation}
and for every locally Lipschitz function
$\varphi\colon \overline{U}\setminus G\to{\mathbb{R}}$ the following 
property
holds:
\begin{equation}
\label{eqAUfp}
\int\limits_{\partial' U}\varphi\,\langle \A(x,\nabla f)\,,
\overline{\bf n}\rangle
\,d{\mathcal H}^{n-1}=
\int\limits_{U}\langle \A(x,\nabla f)\,,\nabla\varphi\,\rangle\,d{\mathcal H}^n\,.
\end{equation}
Here $\overline{\bf n}$ is the unit normal vector of
$\partial' U$
and $d{\mathcal H}^n$ is the volume element on ${\mathbb{R}}^n$.
\end{defn}

\begin{defn}
\label{def2}
We say that a locally Lipschitz  function $f\colon D\to{\mathbb{R}}$ is 
a generalized solution of (\ref{eq2.19}) with the boundary condition
\begin{equation}
\label{areeq3}
f\,\langle \A(x,\nabla f), \overline{\bf n}\rangle =0\,,
\quad x\in
\partial D\setminus G\,,
\end{equation}
if for every subdomain $U\in (G,D)$ with (\ref{eqphiU}),
and for every locally Lipschitz function
$\varphi\colon \overline{U}\setminus G\to{\mathbb{R}}$ the following 
property holds:
\begin{equation}
\label{eeqAUfp}
\int\limits_{\partial' U}\varphi\,f\,\langle \A(x,\nabla f)\,,
\overline{\bf n}\rangle
\,d{\mathcal H}^{n-1}=
\int\limits_{U}\langle \A(x,\nabla f)\,,\nabla(\varphi\,f)\,\rangle\,d{\mathcal H}^n\,.
\end{equation}
\end{defn}

In the case of a smooth boundary  $\partial D$,
and $f\in C^2(D)$, the relation (\ref{eqAUfp}) implies
(\ref{eq2.19}) with (\ref{areq3}) everywhere on $\partial D
\setminus G$. This requirement (\ref{eeqAUfp}) implies
(\ref{eq2.19}) with (\ref{areeq3}) on $\partial D \setminus G$.
See \cite[Section 9.2.1]{[3]}.

The surface integrals exist by (\ref{eqphiU}). Indeed, this
assumption guarantees that $\nabla f(x)$ exists ${\mathcal
H}^{n-1}$-a.e. on $\partial' U$. The assumption that $U\in (G,D)$
implies existence of a normal vector $\overline{\bf n}$ for
${\mathcal H}^{n-1}$-a.e. points on $\partial' U$  \cite[Chapter 2
Section 3.2]{Fe}. Thus,  the scalar product  $\langle \A(x,\nabla
f)\,,\overline{\bf n}\rangle$ is defined and finite a.e. on
$\partial' U$.

\section{Saint-Venant principle}{}

In this section, we will prove the Saint-Venant principle for
solutions of the $\A$-Laplace equation. The Saint-Venant principle states that strains in a body produced by application of a force onto a small part of its surface are of negligible magnitude at distances that are large compared to the diameter of the part where the force is applied. This well known result in elasticity theory is often stated and used in a loose form. For mathematical investigation of the results of this type, see e.g. \cite{BT}.

In this paper the inequalities of the form (\ref{eqIlI}), (\ref{eeqIlI}) are called the Saint-Venant principle (see also, \cite{OI,BT}). Here we consider only the case of canonical domains. We plan to consider the general case in another article.

Let $0<k<n$. Fix a domain $D_0$ in ${\mathbb R}^k$ with compact and smooth boundary,
and write
\[
D=D_0 \times \R^{n-k} = \{x \in {\mathbb R}^n : (x_1,\ldots,x_{k}) \in D_0\}.
\]
We write ${\mathcal P}=\{x\in \partial D:p_k(x)=\alpha\}$ and ${\mathcal Q}=\{x\in \partial
D:p_k(x)=\beta\}$ and $G={\mathcal P}\cup {\mathcal Q}$.
Let $t,\tau\in (\alpha,\beta)$, $t<\tau$, and
$$
\Delta^k(t,\tau)=\{x\in D : t<p_k(x)<\tau\}\,.
$$
For $s\geq 0$ we set
$$
\sigma^k(s)=\big\{x\in \Delta^k(0,\infty):p_k(x)=s\big\}\,.
$$

\begin{thm}
Let $\alpha<\tau'<\tau''<\beta$, and let $0<k<n$.
If $f\colon D\to \R$ is a generalized solution  of {\rm (\ref{eq2.19})} with
the generalized boundary condition {\rm (\ref{areq3})}
on $\partial D\setminus G$, then the inequality 
\begin{equation}
\label{eeqIlI}
I_1(t,\tau')+C_1(t)/\nu_1\le \big(I_1(t,\tau'')+C_1(t)/\nu_1\big)\,
\exp\bigg[-{\frac{\nu_1}{\nu_2}}\,\int\limits_{\tau'}^{\tau''}
\mu_p(\sigma^k(\tau))\,d\tau\bigg]\,.
\end{equation}
holds for all $t\in(\alpha,\tau']$.

If $f\colon D\to \R$ is a generalized solution of {\rm (\ref{eq2.19})} 
with the generalized boundary condition {\rm (\ref{areeq3})} then
\begin{equation}
\label{eqIlI}
I_1(t,\tau')+C_2(t)/\nu_1\le \big(I_1(t,\tau'')+C_2(t)/\nu_1\big)\,
\exp\bigg[-{\frac{\nu_1}{\nu_2}}\,\int\limits_{\tau'}^{\tau''}
\lambda_{p,Z_f(\tau)}^{1/p}(\sigma^k(\tau))\,d\tau\bigg]\,
\end{equation}
holds for all $t\in(\alpha,\tau']$.
Here
$$
I_1(t,\tau)=\int\limits_{\Delta^k(t,\tau)}|\nabla f|^p\,d{\mathcal
H}^{n}\,,
$$
and
\begin{equation}
\label{Zf}
Z_f(\tau)=\{x \in \Sigma_k^*(\tau)\cap \partial D : \lim_{y\to x} 
f(y)=0\}.
\end{equation}
\end{thm}

\proof
{\it The  case A.}
At first we consider  the case in which $f$ is a generalized solution
of (\ref{eq2.19}) with the generalized boundary condition (\ref{areeq3})
on $\partial D\setminus G$.
It is easy to see that a.e. on $D^k_{\alpha,\beta}$,
$$
|\nabla p_k(x)|=1\,.
$$
The domain $\Delta^k(t,\tau)$ belongs to $(G,D)$. Let
$\varphi\colon \overline{U}\setminus G \to \R$ be a locally Lipschitz
function. By (\ref{eeqAUfp}) we have
$$
\int\limits_{\partial' \Delta^k(t,\tau)}\varphi\,f\,\langle \A(x,\nabla f)\,,
\overline{\bf n}\rangle
\,d{\mathcal H}^{n-1}=
\int\limits_{\Delta^k(t,\tau)}\langle \A(x,\nabla f)\,,\nabla(\varphi\,f)\,
\rangle\,d{\mathcal H}^n\,.
$$
But
$$
\partial' \Delta^k(t,\tau)=\sigma^k(t)\cup\sigma^k(\tau)\,.
$$
For $\varphi\equiv 1$, we have by (\ref{eq2.23}) and (\ref{eeqAUfp})
\begin{eqnarray*}
\nu_1\,I_1(t,\tau) &\le& \int\limits_{\Delta^k(t,\tau)}\langle
\A(x,\nabla
f),
\nabla f\rangle\,d{\mathcal H}^n\\
&=&\int\limits_{\sigma^k(\tau)}f\,\langle \A(x,\nabla
f),\nabla p_k(x)\rangle\,d{\mathcal H}^{n-1}-
\int\limits_{\sigma^k(t)} f\,\langle \A(x,\nabla f),\nabla
p_k(x)\rangle\,d{\mathcal H}^{n-1}\,,
\end{eqnarray*}
since $\overline{\bf n}=\nabla p_k(x)\,$ for $x\in\sigma^k(\tau)$
and $\overline{\bf n}=-\nabla p_k(x)\,$ for $x\in\sigma^k(t)$.
We obtain 
\begin{equation}
\label{eqIbc}
\nu_1\,I_1(t,\tau)+C_2(t)\le
\int\limits_{\sigma^k(\tau)}f\,\langle \A(x,\nabla f),\nabla
p_k(x)\rangle\,d{\mathcal H}^{n-1}
\end{equation}
where
$$
C_2(t)=\int\limits_{\sigma^k(t)} f\,\langle \A(x,\nabla f),\nabla
p_k(x)\rangle\,
\,d{\mathcal H}^{n-1}\,.
$$
Note that we may also choose
\begin{equation}
\label{eqIbc_alt}
\tilde{C}_2(\tau)=-\int\limits_{\sigma^k(\tau)} f\,\langle \A(x,\nabla
f),\nabla p_k(x)\rangle\,
\,d{\mathcal H}^{n-1}\,,
\end{equation}
to obtain an inequality similar to (\ref{eqIbc}).

Next we will estimate the right side of (\ref{eqIbc}). By
(\ref{eq2.23}) and the H\"older inequality
\begin{multline*}
\bigg|\int\limits_{\sigma^k(\tau)}f\,\langle \A(x,\nabla f),\nabla
p_k(x)\rangle\,d{\mathcal H}^{n-1}\bigg|\\
\le \int\limits_{\sigma^k(\tau)}|f|\,|\A(x,\nabla f)|\,d{\mathcal
H}^{n-1}
\le \nu_2 \int\limits_{\sigma^k(\tau)}|f|\,|\nabla
f|^{p-1}\,d{\mathcal H}^{n-1}\\
\le
\nu_2\bigg(\int\limits_{\sigma^k(\tau)}|f|^p\,d{\mathcal
H}^{n-1}\bigg)^{1/p}
\bigg(\int\limits_{\sigma^k(\tau)}|\nabla f|^p\,d{\mathcal
H}^{n-1}\bigg)^{(p-1)/p}\,.
\end{multline*}
By using (\ref{eqla1}) we may write
\begin{equation}
\int\limits_{\sigma^k(\tau)}|f|^p\,d{\mathcal H}^{n-1}\le
\lambda_{p,Z_f(\tau)}^{-1}(\sigma^k(\tau))
\,\int\limits_{\sigma^k(\tau)}|\nabla f|^p\,d{\mathcal H}^{n-1}
\end{equation}
and
$$
\bigg|\int\limits_{\sigma^k(\tau)}f\,\langle \A(x,\nabla f),\nabla
p_k(x)\rangle\,d{\mathcal H}^{n-1}\bigg|\le
\nu_2\,\lambda_{p,Z_f(\tau)}^{-1/p}(\sigma^k(\tau))\,\int\limits_{\sigma^k(\tau)}|\nabla
f|^p\,d{\mathcal H}^{n-1}\,.
$$

By (\ref{eqIbc}) and the Fubini theorem
$$
\nu_1\,I_1(t,\tau)+C_2(t)\le \nu_2\,\lambda_{p,Z_f(\tau)}^{-1/p}(\sigma^k(\tau))\,
{\frac{dI_1}{d\tau}}(t,\tau)
$$
and
$$
{\frac{\nu_1}{\nu_2}}\,\lambda_{p,Z_f(\tau)}^{1/p}(\sigma^k(\tau))\,\le
{\frac{dI_1}{d\tau}}(t,\tau)\bigg/\bigg(I_1(t,\tau)
+{\frac{C_2(t)}{\nu_1}}\bigg)\,.
$$
By integrating this differential inequality we have
$$
\exp\bigg[{\frac{\nu_1}{\nu_2}}\,\int\limits_{\tau'}^{\tau''}
\lambda_{p,Z_f(\tau)}^{1/p}(\sigma^k(\tau))\,d\tau\bigg]\le
{\frac{I_1(t,\tau'')+C_2(t)/\nu_1}{I_1(t,\tau')+C_2(t)/\nu_1}}
$$
for arbitrary $\tau',\tau''\in(\alpha,\beta)$ with $\tau'<\tau''$.
We have shown that
\begin{equation}
I_1(t,\tau')+C_2(t)/\nu_1\le \big(I_1(t,\tau'')+C_2(t)/\nu_1\big)\,
\exp\bigg[-{\frac{\nu_1}{\nu_2}}\,\int\limits_{\tau'}^{\tau''}
\lambda_{p,Z_f(\tau)}^{1/p}(\sigma_p^k(\tau))\,d\tau\bigg]\,.
\end{equation}

{\it The  case B.} Now we assume that $f$ is a generalized solution  of
(\ref{eq2.19}) with the boundary condition (\ref{areq3})
on $\partial D\setminus G$. Fix $t<\tau$.
By choosing $\varphi\equiv 1$ in (\ref{eqAUfp}) we
see that
$$
\int\limits_{\sigma^k(t)\cup\sigma^k(\tau)}\langle \A(x,\nabla f)\,,
\overline{\bf n}\rangle
\,d{\mathcal H}^{n-1}=0\,.
$$
For an arbitrary constant $C$, we get from this and (\ref{eqAUfp})
\begin{equation}
\label{eqphfC}
\int\limits_{\sigma^k(t)\cup\sigma^k(\tau)}(f-C)\,\langle \A(x,\nabla f)\,,
\overline{\bf n}\rangle
\,d{\mathcal H}^{n-1}=
\int\limits_{\Delta^k(t,\tau)}\langle \A(x,\nabla f)\,,\nabla f\,\rangle\,d{\mathcal H}^n\,.
\end{equation}
Thus
$$
\int\limits_{\Delta^k(t,\tau)}\langle \A(x,\nabla f)\,,\nabla f\,\rangle\,d{\mathcal H}^n\le
C_1(t)+\int\limits_{\sigma^k(\tau)}|f-C|\,|\A(x,\nabla p_k(x))|\,
d{\mathcal H}^{n-1},
$$
where
$$
C_1(t)=\int_{\sigma^k(t)} |f-C|\,|\A(x,\nabla p_k(x))|\,
d{\mathcal H}^{n-1},
$$
or
\begin{equation}
\label{eeqIbc}
\nu_1\,I_1(t,\tau)+C_1(t)\le
\nu_2\,\int\limits_{\sigma^k(\tau)}|f-C|\,|\nabla f|^{p-1}\,
d{\mathcal H}^{n-1}\,.
\end{equation}
As above, we obtain
\begin{equation}
\label{arbit}
\int\limits_{\sigma^k(\tau)}|f-C|\,|\nabla f|^{p-1}\,
d{\mathcal H}^{n-1}\le \bigg(\int\limits_{\sigma^k(\tau)}|f-C|^p\,
d{\mathcal H}^{n-1}\bigg)^{1/p}\bigg(\int\limits_{\sigma^k(\tau)}
|\nabla f|^p\,
d{\mathcal H}^{n-1}\bigg)^{(p-1)/p}\,.
\end{equation}
By using (\ref{eqla2}) we get
\begin{equation}
\label{c3}
\bigg(\int\limits_{\sigma^k(\tau)}|f-C_3|^p\,
d{\mathcal H}^{n-1}\bigg)^{1/p}\le
\mu^{-1/p}_{p}(\sigma^k(\tau))\,
\bigg(\int\limits_{\sigma^k(\tau)}
|\nabla f|^p\,
d{\mathcal H}^{n-1}\bigg)^{1/p},
\end{equation}
where $C_3=C_3(f)$ is the constant from (\ref{eqla2}). 
Then by (\ref{arbit}) and (\ref{c3}),
$$
\int\limits_{\sigma^k(\tau)}|f-C_3|\,|\nabla f|^{p-1}\,
d{\mathcal H}^{n-1}\le \mu^{-1}_{p}(\sigma^k(\tau))
\int\limits_{\sigma^k(\tau)}
|\nabla f|^p\,d{\mathcal H}^{n-1}\,,
$$
and by (\ref{eeqIbc}) we have
$$
\nu_1\,I_1(t,\tau)+C_1(t)\le \nu_2\,\mu^{-1}_{p}(\sigma^k(\tau))
\int\limits_{\sigma^k(\tau)}
|\nabla f|^p\,d{\mathcal H}^{n-1}
$$
or
$$
\nu_1\,I_1(t,\tau)+C_1(t)\le
\nu_2\,\mu^{-1}_{p}(\sigma^k(\tau))
{\frac{dI_1}{dt}}
(t,\tau)\,.
$$
By integrating this inequality we have shown that
\begin{equation}
I_1(t,\tau')+C_1(t)/\nu_1\le \left(I_1(t,\tau'')+C_1(t)/\nu_1\right)\,
\exp\bigg[-{\frac{\nu_1}{\nu_2}}\,\int\limits_{\tau'}^{\tau''}
\mu_p(\sigma^k(\tau))\,d\tau\bigg]\,.
\end{equation}

\qed

\section{Stagnation zones}{}

Next we apply the Saint-Venant principle to obtain information
about stagnation zones of generalized solutions of the
equation (\ref{eq2.19}).  We first consider zones with respect to the
Sobolev norm. Other results of this type follow immediately from
well-known imbedding theorems.

\subsection*{Stagnation zones with respect to the $W^1_p$-norm}{}

We rewrite (\ref{eeqIlI}) and (\ref{eqIlI}) in another form. 
Let $0<k<n$ and let $0<\alpha<\beta$. Fix a domain $D_0$ in ${\mathbb R}^k$ with compact and smooth boundary, and write
\[
D=D_0 \times \R^{n-k} = \{x \in {\mathbb R}^n : (x_1,\ldots,x_{k}) \in D_0\}.
\]
We write
$$
p^*_k(x)=p_k(x)-{\frac{\alpha+\beta}{2}}\,.
$$
For $x\in D^k_{\alpha,\beta}$ and
\begin{equation}
\label{betastar}
\beta^*={\frac{\beta-\alpha}{2}}
\end{equation}
we have
$$
-\beta^*<p^*_k(x)<\beta^*\,,
$$
and we denote
$$
D^{*,k}_{\beta^*}=\{x\in {\bf R}^n: -\beta^* <p^*_k(x)<\beta^*\}\,.
$$
Let  $-\beta^* <\tau' \le \tau''<\beta^*$. We write
\begin{equation}
\label{deltastar}
\Delta^{*,k}(\tau',\tau'')=\{x\in D :\tau'<p^*_k(x)<\tau''\}
\end{equation}
and
$$
I_2(\tau',\tau'')=\int\limits_{\Delta^{*,k}(\tau',\tau'')}|\nabla
f|^p\, d{\mathcal H}^n\,.
$$

Let $0<\tau'<\tau''<\beta^*.$ By (\ref{eqIlI}) we have for
$t\in(-\tau,\tau)$
$$
I_2(t,\tau')+C_4(t)/\nu_1\le \big(I_2(t,\tau'')+C_4(t)/\nu_1\big)\,
\exp\bigg[-{\frac{\nu_1}{\nu_2}}\,\int\limits_{\tau'}^{\tau''}
\lambda_{p,Z^*_f(\tau)}^{1/p}(\sigma^{*,k}(\tau))\,d\tau\bigg]\,,
$$
where
\begin{equation}
\label{Zfstar}
Z^*_f(\tau)=\{x \in \partial D : p_k^*(x)=\tau \land \lim_{y\to x} 
f(y)=0\}.
\end{equation}

By choosing the estimate as in (\ref{eqIbc_alt}) we also have
$$
I_2(-\tau',t)+\tilde{C}_4(t)/\nu_1\le
\big(I_2(-\tau'',t)+\tilde{C}_4(t)/\nu_1\big)\,
\exp\bigg[-{\frac{\nu_1}{\nu_2}}\,\int\limits_{-\tau''}^{-\tau'}
\lambda_{p,Z^*_f(\tau)}^{1/p}(\sigma^{*,k}(\tau))\,d\tau\bigg]\,,
$$
where
\begin{equation}
\label{sigmastar}
\sigma^{*,k}(s)=\{x\in \Delta^{*,k}(-\infty,\infty):p^*_k(x)=s\}\,.
\end{equation}
By adding these inequalities and noting that $C_4(t)+\tilde{C}_4(t)=0$
we
obtain
\begin{multline*}
I_2(-\tau',t)+ I_2(t,\tau')\le
\big(I_2(-\tau'',t)+I_2(t,\tau'')\big)\\
\times
\max\Bigg\{\exp\bigg[-{\frac{\nu_1}{\nu_2}}\,\int\limits_{-\tau''}^{-\tau'}
\lambda_{p,Z^*_f(\tau)}^{1/p}(\sigma^{*,k}(\tau))\,d\tau\bigg],
\exp\bigg[-{\frac{\nu_1}{\nu_2}}\,\int\limits_{\tau'}^{\tau''}
\lambda_{p,Z^*_f(\tau)}^{1/p}(\sigma^{*,k}(\tau))\,d\tau\bigg]
\Bigg\}\,.
\end{multline*}
Thus we have the estimate
\begin{multline}
\label{eqSVP1}
I_2(-\tau',\tau')\le I_2(-\tau'',\tau'')\\
\times
\max\Bigg\{\exp\bigg[-\frac{\nu_1}{\nu_2}\,\int\limits_{-\tau''}^{-\tau'}
\lambda_{p,Z^*_f(\tau)}^{1/p}(\sigma^{*,k}(\tau))\,d\tau\bigg],
\exp\bigg[-\frac{\nu_1}{\nu_2}\,\int\limits_{\tau'}^{\tau''}
\lambda_{p,Z^*_f(\tau)}^{1/p}(\sigma^{*,k}(\tau))\,d\tau\bigg]
\Bigg\}\,.
\end{multline}
Similarly, from (\ref{eeqIlI}) we obtain
\begin{multline}
\label{eqSVP2}
I_2(-\tau',\tau')\le I_2(-\tau'',\tau'')\\
\times
\max\Bigg\{\exp\bigg[-{\frac{\nu_1}{\nu_2}}\,\int\limits_{-\tau''}^{-\tau'}
\mu_p(\sigma^{*,k}(\tau))\,d\tau\bigg],
\exp\bigg[-{\frac{\nu_1}{\nu_2}}\,\int\limits_{\tau'}^{\tau''}
\mu_p(\sigma^{*,k}(\tau))\,d\tau\bigg]
\Bigg\}\,.
\end{multline}

From this we obtain the following theorem on stagnation $W^1_p$-zones:

\begin{thm}
Let $0<k<n$, $\beta>\alpha>0$, and let $-\beta^* <\tau' \le
\tau''<\beta^*$ where $\beta^*$ is as in {\rm (\ref{betastar})}. If $f$
is a solution of {\rm (\ref{eq2.19})} on $D$
 with the generalized  boundary condition {\rm (\ref{areq3})}
on $\partial D\setminus G$, where
$G=\{x\in \partial D:p^*_k(x)=\pm\beta^*\}$
and
$$
\max\Bigg\{\exp\bigg[-{\frac{\nu_1}{\nu_2}}\,\int\limits_{-\tau''}^{-\tau'}
\mu_p(\sigma^{*,k}(\tau))\,d\tau\bigg],
\exp\bigg[-{\frac{\nu_1}{\nu_2}}\,\int\limits_{\tau'}^{\tau''}
\mu_p(\sigma^{*,k}(\tau))\,d\tau\bigg]
\Bigg\}<s^{1/p}\,,
$$
or a solution of
{\rm (\ref{eq2.19})} on 
$D$ with the generalized
boundary condition {\rm (\ref{areeq3})} on $\partial D\setminus G$ and
$$
\max\Bigg\{\exp\bigg[-{\frac{\nu_1}{\nu_2}}\,\int\limits_{-\tau''}^{-\tau'}
\lambda_{p,Z^*_f(\tau)}^{1/p}(\sigma^{*,k}(\tau))\,d\tau\bigg],
\exp\bigg[-{\frac{\nu_1}{\nu_2}}\,\int\limits_{\tau'}^{\tau''}
\lambda_{p,Z^*_f(\tau)}^{1/p}(\sigma^{*,k}(\tau))\,d\tau\bigg]
\Bigg\} <s^{1/p}\,,
$$
then the subdomain $\Delta^{*,k}(-\tau',\tau')$ is a $s$-zone with
respect to the $W^1_p$-norm i.e.,
$$
\int\limits_{\Delta^{*,k}(-\tau',\tau')}|\nabla f|^p\,d{\mathcal
H}^n<s\,,
$$
where $\Delta^{*,k}$ is as in {\rm (\ref{deltastar})}.
\end{thm}


\subsection*{Stagnation zones with respect to the $L^p$-norm}{}
Let $\beta>\alpha>0$, and let $-\beta^* <\tau' \le \tau''<\beta^*$ where $\beta^*$ is as in (\ref{betastar}).

Denote by $C_5$ the best constant of the imbedding
theorem from $W^1_p(D^{*,k}_{\beta^*})$ to $L^p(D^{*,k}_{\beta^*})$, i.e. in the inequality
$$
\|g-C\|_{L^p(D^{*,k}_{\beta^*})}\le
C_5\,\|g\|_{W^1_p(D^{*,k}_{\beta^*})}\,,
$$
if such constant exist (see Maz'ya \cite{Mazya:1985} or \cite{Adams}). Then
 we obtain from (\ref{eqSVP1}), (\ref{eqSVP2})
\begin{multline}
\label{eqSVP3}
\|f-C\|^p_{L^p(\Delta^{*,k}(-\tau',\tau'))}\le
C_5^{p}\,I_2(-\tau'',\tau'')\\
\times
\max\Bigg\{\exp\bigg[-{\frac{\nu_1}{\nu_2}}\,\int\limits_{-\tau''}^{-\tau'}
\lambda_{p,Z^*_f(\tau)}^{1/p}(\sigma^{*,k}(\tau))\,d\tau\bigg],
\exp\bigg[-{\frac{\nu_1}{\nu_2}}\,\int\limits_{\tau'}^{\tau''}
\lambda_{p,Z^*_f(\tau)}^{1/p}(\sigma^{*,k}(\tau))\,d\tau\bigg] \Bigg\}\,,
\end{multline}
and
\begin{multline}
\label{eqSVP4}
\|f-C\|^p_{L^p(\Delta^{*,k}_{\tau'})}\le
C_5^{p} I_2(-\tau'',\tau'')\\
\times
\max\Bigg\{\exp\bigg[-{\frac{\nu_1}{\nu_2}}\,\int\limits_{-\tau''}^{-\tau'}
\mu_p(\sigma^{*,k}(\tau))\,d\tau\bigg],
\exp\bigg[-{\frac{\nu_1}{\nu_2}}\,\int\limits_{\tau'}^{\tau''}
\mu_p(\sigma^{*,k}(\tau))\,d\tau\bigg] \Bigg\}\,.
\end{multline}
These relations can be used to obtain information about
stagnation zones with respect to the $L^p$-norm. Namely, we have:

\begin{thm}
Let $0<k<n$, and let
\[
D=D_0 \times \R^{n-k} = \{x \in {\mathbb R}^n : (x_1,\ldots,x_{k}) \in D_0\},
\]
where $D_0$ is a domain in ${\mathbb R}^k$ with compact and smooth boundary.
If $f$ is a solution of {\rm (\ref{eq2.19})} on $D$, with the generalized boundary
condition {\rm (\ref{areq3})} (or {\rm (\ref{areeq3})}) on $\partial D\setminus G$, where $G=\{x\in \partial D:p^*_k(x)=\pm\beta^*\}$, and the right side of
{\rm (\ref{eqSVP3})} (or {\rm (\ref{eqSVP4}))} is smaller
than $s>0$, then the domain $\Delta^{*,k}(-\tau',\tau')$ is a stagnation
zone  with the deviation $s^p$ in the sense of the $L^p$-norm on
$D$.
\end{thm}


\subsection*{Stagnation zones for bounded, uniformly continuous
functions}
Let $\beta>\alpha>0$, and let $-\beta^* <\tau' \le \tau''<\beta^*$ where $\beta^*$ is as in (\ref{betastar}).

As before, denote by $C_6$ the best constant of the
imbedding theorem from $W^1_p(D^{*,k}_{\beta^*})$ to $C(D^{*,k}_{\beta^*})$, i.e. in the inequality
\begin{equation}
\label{sob2}
\|g-C\|_{C(D^{*,k}_{\beta^*})}\le
C_6\,\|g\|_{W^1_p(D^{*,k}_{\beta^*})}\,,
\end{equation}
if such constant exists. For example, if the domain $D^{*,k}_{\beta^*}$ is convex, then (\ref{sob2}) holds for $p>n$ (see Maz'ya \cite{Mazya:1985} or \cite[p.85]{Adams}).

In this case from (\ref{eqSVP1}), (\ref{eqSVP2}) we obtain
\begin{multline}
\label{eqSVP5}
\|f-C\|_{C(\Delta^{*,k}(-\tau',\tau'))}\le
C_6^{p}\,I_2(-\tau'',\tau'')\\
\times
\max\Bigg\{\exp\bigg[-{\frac{\nu_1}{\nu_2}}\,\int\limits_{-\tau''}^{-\tau'}
\lambda_{p,Z^*_f(\tau)}^{1/p}(\sigma^{*,k}(\tau))\,d\tau\bigg],
\exp\bigg[-{\frac{\nu_1}{\nu_2}}\,\int\limits_{\tau'}^{\tau''}
\lambda_{p,Z^*_f(\tau)}^{1/p}(\sigma^{*,k}(\tau))\,d\tau\bigg]
\Bigg\}\,,
\end{multline}
and
\begin{multline}
\label{eqSVP6}
\|f-C\|_{C(\Delta^{*,k}(-\tau',\tau'))}\le
C_6^{p}I_2(-\tau'',\tau'')\\
\times
\max\Bigg\{\exp\bigg[-{\frac{\nu_1}{\nu_2}}\,\int\limits_{-\tau''}^{-\tau'}
\mu_p(\sigma^{*,k}(\tau))\,d\tau\bigg],
\exp\bigg[-{\frac{\nu_1}{\nu_2}}\,\int\limits_{\tau'}^{\tau''}
\mu_p(\sigma^{*,k}(\tau))\,d\tau\bigg]
\Bigg\}\,.
\end{multline}
These relations can be used to obtain theorems about
stagnation zones for bounded uniformly continuous functions. 

\begin{thm}
Let $0<k<n$. If $f$ is a solution of
{\rm (\ref{eq2.19})}, $p>n$, on $D$ where $D$ is as before with
the generalized boundary condition {\rm (\ref{areq3})} (or {\rm
(\ref{areeq3})})
on $\partial D\setminus G$  where
$G=\{x\in \partial D:p^*_k(x)=\pm\beta^*\}$, and the right side of
{\rm (\ref{eqSVP5})} (or
{\rm (\ref{eqSVP6}))} is smaller
than $s>0$, then the domain $\Delta^{*,k}(-\tau',\tau')$ is a
stagnation zone with the deviation $s$ in the sense of the norm
$\|\cdot\|_{C^0(\Delta^{*,k}(-\tau',\tau'))}$.
\end{thm}


\section{Other applications}{}

Next we prove Phragm\'en-Lindel\"of type theorems for the solutions
of the $\A$-Laplace equation with
boundary conditions (\ref{areq3}) and (\ref{areeq3}).

\subsection*{Estimates for $W^1_p$-norms}{}
Let $\beta>\alpha>0$, and let $D_0$ be a domain in ${\mathbb R}^k$ with compact and smooth boundary. Write
\[
D=D_0 \times \R^{n-k} = \{x \in {\mathbb R}^n : (x_1,\ldots,x_{k}) \in D_0\}.
\]
Suppose that $\beta^*$ is as in (\ref{betastar}).
First we will prove some estimates of the $W^1_p$-norm of a solution.
Let $f$ be a solution of
{\rm (\ref{eq2.19})} on $D^{*,k}_{\beta^*}$ with the generalized
boundary condition
{\rm (\ref{areq3})} on $\partial D\setminus G$. Fix
$0<\tau'<\tau''<\beta^*$ and estimate
$\|f\|_{W^1_p(\Delta^{*,k}(-\tau',\tau'))}$.

Let $\psi:[\tau',\tau'']\to (0,\infty)$ be a Lipschitz function such
that
\begin{equation}
\label{psi}
\psi(\tau')=1\,,\quad \psi(\tau'')=0\,.
\end{equation}
We choose
\begin{equation}
\label{eq116}
\phi(t)=\left\{
\begin{array}{lll}
1&\text{ for } |t|<\tau',\\
\psi(|t|)&\text{ for } \tau'\le |t| \le \tau''.\\
\end{array}
\right.
\end{equation}
The function $\varphi(x)=\phi(p^*_k(x))$ is admissible in Definition
\ref{def1}
for
$$
U=\Delta^{*,k}(-\tau'',\tau'')\,.
$$
As in (\ref{eqphfC}) we may by (\ref{eqAUfp}) write
\begin{multline*}
\int\limits_{\sigma^{*,k}(-\tau'')\cup\sigma^{*,k}(\tau'')}\phi^p(p^*_k(x))(f-C)\,\langle
\A(x,\nabla f)\,,
\overline{\bf n}\rangle
\,d{\mathcal H}^{n-1}\\
=
\int\limits_{\Delta^{*,k}(-\tau'',\tau'')}\langle
\A(x,\nabla
f)\,,\nabla
\big(\phi^p(p^*_k(x))(f-C)\big)\,\rangle\,d{\mathcal H}^n\,.
\end{multline*}
By the construction of $\phi$, (\ref{psi}) and (\ref{eq116}),
the surface integral is equal to zero, and we have
\begin{multline*}
\int\limits_{\Delta^{*,k}(-\tau'',\tau'')}\phi^p(p^*_k(x))\langle
\A(x,\nabla
f)\,,\nabla f\rangle\,d{\mathcal H}^n\\
=-p\int\limits_{\Delta^{*,k}(-\tau'',\tau'')}
\phi^{p-1}(p^*_k(x))\,(f-C)\,\langle \A(x,\nabla f)\,,
\nabla \phi(p^*_k(x))\,\rangle\,d{\mathcal H}^n\,.
\end{multline*}
Thus by (\ref{eq2.23})
\begin{multline}
\label{eq2239}
\nu_1 \int\limits_{\Delta^{*,k}(-\tau'',\tau'')}\phi^p(p^*_k(x))|\nabla
f|^p\,d{\mathcal H}^n\\
\le
p\nu_2\int\limits_{\Delta^{*,k}(-\tau'',\tau'')}
\phi^{p-1}(p^*_k(x))\,|f-C|\,|\nabla f|^{p-1}\,
|\nabla \phi(p^*_k(x))|\,d{\mathcal H}^n\,.
\end{multline}
Now we note that
$$
|\nabla \phi(p^*_k(x))|=|\phi'(p^*_k(x))|
$$
and by the H\"older inequality
\begin{multline*}
\int\limits_{\Delta^{*,k}(-\tau'',\tau'')}
\phi^{p-1}(p^*_k(x))\,|f-C|\,|\nabla f|^{p-1}\,
|\nabla \phi(p^*_k(x))|\,d{\mathcal H}^n\\
=\int\limits_{\Delta^{*,k}(-\tau'',\tau'')}
\phi^{p-1}(p^*_k(x))\,|\nabla f|^{p-1}\,|\phi'(p^*_k(x))|\,|f-C|\,
d{\mathcal H}^n\\
\le
\bigg(\int\limits_{\Delta^{*,k}(-\tau'',\tau'')}
\phi^p(p^*_k(x))\,|\nabla f|^{p}\,
d{\mathcal
H}^n\bigg)^{(p-1)/p}\bigg(\int\limits_{\Delta^{*,k}(-\tau'',\tau'')}
|\phi'(p^*_k(x))|^p\,|f-C|^p\,
d{\mathcal H}^n\bigg)^{1/p}
\end{multline*}
From this inequality and (\ref{eq2239}) we obtain
$$
\nu_1^p\int\limits_{\Delta^{*,k}(-\tau',\tau')}\phi^p(p^*_k(x))
|\nabla f|^p\,d{\mathcal H}^n\le p^p\nu_2^p
\int\limits_{\Delta^{*,k}(-\tau'',\tau'')}
|\phi'(p^*_k(x))|^p\,|f-C|^p\,
d{\mathcal H}^n\,.
$$
Because $\phi(p^*_k(x))\equiv 1$ on $\Delta^{*,k}(-\tau',\tau')$ we
have the following inequality:
\begin{equation}
\label{eq2230}
\nu_1^p\int\limits_{\Delta^{*,k}(-\tau',\tau')}
|\nabla f|^p\,d{\mathcal H}^n\le p^p\nu_2^p
\int\limits_{\Delta^{*,k}(-\tau'',\tau'')\setminus
\Delta^{*,k}(-\tau',\tau')}
|\psi'(p^*_k(x))|^p\,|f-C|^p\,
d{\mathcal H}^n\,.
\end{equation}

Next we will find
$$
\min_{\psi}\int\limits_{\Delta^{*,k}(-\tau'',\tau'')\setminus
\Delta^{*,k}(-\tau',\tau')}|\psi'(p^*_k(x))|^p\,|f-C|^p\,
d{\mathcal H}^n\,,
$$
where the minimum is taken over all $\psi$ in (\ref{eq116}).
We have
\begin{multline*}
\int\limits_{\Delta^{*,k}(-\tau'',\tau'')\setminus
\Delta^{*,k}(-\tau',\tau')}
|\psi'(p^*_k(x))|^p\,|f-C|^p\,
d{\mathcal H}^n\\
=
\int\limits_{-\tau''}^{-\tau'}|\psi'(\tau)|^p\,d\tau
\int\limits_{\sigma^{*,k}(\tau)}
|f(x)-C|^p\,d{\mathcal H}^{n-1}\\
+\int\limits_{\tau'}^{\tau''}|\psi'(\tau)|^p\,d\tau\int\limits_{\sigma^{*,k}(\tau)}
|f(x)-C|^p\,d{\mathcal H}^{n-1}
\end{multline*}
and
\begin{multline}
\label{eq130}
\min_{\psi}\int\limits_{\Delta^{*,k}(-\tau'',\tau'')\setminus
\Delta^{*,k}(-\tau',\tau')}
|\psi'(p^*_k(x))|^p\,|f(x)-C|^p\,
d{\mathcal H}^n \\
\le
\min_{\psi}\int\limits_{-\tau''}^{-\tau'}|\psi'(\tau)|^p\,d\tau
\int\limits_{\sigma^{*,k}(\tau)}
|f(x)-C|^p\,d{\mathcal H}^{n-1} \\
+\min_{\psi}\int\limits_{\tau'}^{\tau''}|\psi'(\tau)|^p\,d\tau
\int\limits_{\sigma^{*,k}(\tau)}
|f(x)-C|^p\,d{\mathcal H}^{n-1}\equiv A_1+A_2\,.
\end{multline}
Because by the H\"older inequality
\begin{multline*}
1\le \bigg(\int\limits_{\tau'}^{\tau''}|\psi'(\tau)|\,d\tau\bigg)^p\le
\bigg[\int\limits_{\tau'}^{\tau''}|\psi'(\tau)|^pd\tau\,
\int\limits_{\sigma^{*,k}(\tau)}
|f(x)-C|^p\,d{\mathcal H}^{n-1}\bigg]\\
\times
\Bigg[\int\limits_{\tau'}^{\tau''}d\tau\,\bigg(\int\limits_{\sigma^{*,k}(\tau)}
|f(x)-C|^p\,d{\mathcal
H}^{n-1}\bigg)^{1/(1-p)}\Bigg]^{p-1}\,,
\end{multline*}
we have
\begin{multline*}
\Bigg[\int\limits_{\tau'}^{\tau''}d\tau\,\bigg(\int\limits_{\sigma^{*,k}(\tau)}
|f(x)-C|^p\,d{\mathcal H}^{n-1}\bigg)^{1/(1-p)}\Bigg]^{1-p}\\
\le
\int\limits_{\tau'}^{\tau''}|\psi'(\tau)|^pd\tau\,\int\limits_{\sigma^{*,k}(\tau)}
|f(x)-C|^p\,d{\mathcal H}^{n-1},
\end{multline*}
and hence
$$
A_2\ge
\Bigg[\int\limits_{\tau'}^{\tau''}d\tau\,\bigg(\int\limits_{\sigma^{*,k}(\tau)}
|f(x)-C|^p\,d{\mathcal H}^{n-1}\bigg)^{1/(1-p)}\Bigg]^{1-p}\,.
$$
It is easy to see that here the equality holds for a
special choice of $\psi$. Thus
$$
A_2=\Bigg[\int\limits_{\tau'}^{\tau''}d\tau\,\bigg(\int\limits_{\sigma^{*,k}(\tau)}
|f(x)-C|^p\,d{\mathcal H}^{n-1}\bigg)^{1/(1-p)}\Bigg]^{1-p}\,.
$$
Similarly,
$$
A_1=\Bigg[\int\limits_{-\tau''}^{-\tau'}d\tau\,\bigg(\int\limits_{\sigma^{*,k}(\tau)}
|f(x)-C|^p\,d{\mathcal H}^{n-1}\bigg)^{1/(1-p)}\Bigg]^{1-p}\,.
$$
From (\ref{eq130}) we obtain
\begin{multline*}
\min_{\psi}\int\limits_{\Delta^{*,k}(-\tau'',\tau'')\setminus
\Delta^{*,k}(-\tau',\tau')}
|\psi'(p^*_k(x))|^p\,|f-C|^p\,
d{\mathcal H}^n \\
\le\Bigg[\int\limits_{-\tau''}^{-\tau'}d\tau\,
\bigg(\int\limits_{\sigma^{*,k}(\tau)}
|f(x)-C|^p\,d{\mathcal H}^{n-1}\bigg)^{1/(1-p)}\Bigg]^{1-p}\\
+
\Bigg[\int\limits_{\tau'}^{\tau''}d\tau\,\bigg(\int\limits_{\sigma^{*,k}(\tau)}
|f(x)-C|^p\,d{\mathcal H}^{n-1}\bigg)^{1/(1-p)}\Bigg]^{1-p}\,.
\end{multline*}
By using (\ref{eq2230}) we obtain the inequality:
\begin{multline*}
p^{-p}\Big({\frac{\nu_1}{\nu_2}}\Big)^p\int\limits_{\Delta^{*,k}(-\tau',\tau')}
|\nabla f|^p\,d{\mathcal H}^n
\le\Bigg[\int\limits_{-\tau''}^{-\tau'}d\tau\,\bigg(\int\limits_{\sigma^{*,k}(\tau)}
|f(x)-C|^p\,d{\mathcal H}^{n-1}\bigg)^{1/(1-p)}\Bigg]^{1-p} \\
+\Bigg[\int\limits_{\tau'}^{\tau''}d\tau\,\bigg(\int\limits_{\sigma^{*,k}(\tau)}
|f(x)-C|^p\,d{\mathcal H}^{n-1}\bigg)^{1/(1-p)}\Bigg]^{1-p}\,,
\end{multline*}
where $C$ is an arbitrary constant. From this we obtain
\begin{multline}
\label{eq203}
\int\limits_{\Delta^{*,k}(-\tau',\tau')}
|\nabla f|^p\,d{\mathcal H}^n
\le
C_7
\max\Bigg\{\Bigg[\int\limits_{-\tau''}^{-\tau'}d\tau\,\bigg(\int\limits_{\sigma^{*,k}(\tau)}
|f(x)-C|^p\,d{\mathcal H}^{n-1}\bigg)^{1/(1-p)}\Bigg]^{1-p},\\
\Bigg[\int\limits_{\tau'}^{\tau''}d\tau\,\bigg(\int\limits_{\sigma^{*,k}(\tau)}
|f(x)-C|^p\,d{\mathcal H}^{n-1}\bigg)^{1/(1-p)}\Bigg]^{1-p}
\Bigg\}\,
\end{multline}
where $C_7=2p^{p}(\nu_2/\nu_1)^p$.

Similarly, for the solutions of the $\A$-Laplace equation with the boundary
condition
(\ref{areeq3}) we may prove that
\begin{multline*}
p^{-p}\bigg({\frac{\nu_1}{\nu_2}}\bigg)^p\int\limits_{\Delta^{*,k}(-\tau',\tau')}
|\nabla f|^p\,d{\mathcal H}^n
\le\Bigg[\int\limits_{-\tau''}^{-\tau'}d\tau\,\bigg(\int\limits_{\sigma^{*,k}(\tau)}
|f(x)|^p\,d{\mathcal H}^{n-1}\bigg)^{1/(1-p)}\Bigg]^{1-p}\\
+\Bigg[\int\limits_{\tau'}^{\tau''}d\tau\,\bigg(\int\limits_{\sigma^{*,k}(\tau)}
|f(x)|^p\,d{\mathcal H}^{n-1}\bigg)^{1/(1-p)}\Bigg]^{1-p}\,.
\end{multline*}
It follows that
\begin{multline}
\label{eeq203}
\int\limits_{\Delta^{*,k}(-\tau',\tau')}
|\nabla f|^p\,d{\mathcal H}^n
\le C_7
\max\Bigg\{\Bigg[\int\limits_{-\tau''}^{-\tau'}d\tau\,\bigg(\int\limits_{\sigma^{*,k}(\tau)}
|f(x)|^p\,d{\mathcal H}^{n-1}\bigg)^{1/(1-p)}\Bigg]^{1-p},\\
\Bigg[\int\limits_{\tau'}^{\tau''}d\tau\,\bigg(\int\limits_{\sigma^{*,k}(\tau)}
|f(x)|^p\,d{\mathcal H}^{n-1}\bigg)^{1/(1-p)}\Bigg]^{1-p}
\Bigg\}\,.
\end{multline}

\subsection*{Phragm\'en-Lindel\"of type theorems I}

We prove Phragm\'en-Lindel\"of type theorems for cylindrical domains.
Let $k=n-1$. Fix a domain $D_0$ in ${\mathbb R}^{n-1}$ with compact and smooth boundary. Consider the domain
\[
D=D_0 \times \R = \{x \in {\mathbb R}^n : (x_1,\ldots,x_{n-1}) \in D_0\}.
\]
Let $f_0\colon D\to \R$ be a generalized solution of (\ref{eq2.19}) with
(\ref{areq5}) and (\ref{eq2.23}) satisfying the boundary condition
(\ref{areq3}) on $\partial D$.

Fix $\beta>\alpha>0$, and let $\beta^*$ be as in (\ref{betastar}). Let
$f(x)=f_0(x-\beta^* e_n)$, where $e_n$ is the $n$:th unit coordinate
vector, and let $0<\tau'<\tau''<\beta^*<\infty\,.$
By (\ref{eq203})
\begin{multline*}
\int\limits_{\Delta^{*,k}(-\tau'',\tau'')}
|\nabla f|^p\,d{\mathcal H}^n
\le C_7 
\max\Bigg\{\Bigg[\int\limits_{-\tau''-1}^{-\tau''}d\tau\,
\bigg(\int\limits_{\sigma^{*,n-1}(\tau)}
|f(x)-C|^p\,d{\mathcal H}^{n-1}\bigg)^{1/(1-p)}\Bigg]^{1-p},
\\
\Bigg[\int\limits_{\tau''}^{\tau''+1}d\tau\,
\bigg(\int\limits_{\sigma^{*,n-1}(\tau)}
|f(x)-C|^p\,d{\mathcal H}^{n-1}\bigg)^{1/(1-p)}\Bigg]^{1-p}
\Bigg\}\,.
\end{multline*}
By using (\ref{eqSVP2}) we obtain from this the inequality:
\begin{multline*}
I_2(-\tau',\tau')\le C_7 
\max\Bigg\{\Bigg[\int\limits_{-\tau''-1}^{-\tau''}d\tau\,
\bigg(\int\limits_{\sigma^{*,n-1}(\tau)}
|f(x)-C|^p\,d{\mathcal H}^{n-1}\bigg)^{1/(1-p)}\Bigg]^{1-p},
\\
\Bigg[\int\limits_{\tau''}^{\tau''+1}d\tau\,\bigg(
\int\limits_{\sigma^{*,n-1}(\tau)}
|f(x)-C|^p\,d{\mathcal H}^{n-1}\bigg)^{1/(1-p)}\Bigg]^{1-p}
\Bigg\}\,\\
\times
\max\Bigg\{\exp\bigg[-{\frac{\nu_1}{\nu_2}}\,\int\limits_{-\tau''}^{-\tau'}
\mu_p(\sigma^{*,n-1}(\tau))\,d\tau\bigg],
\exp\bigg[-{\frac{\nu_1}{\nu_2}}\,\int\limits_{\tau'}^{\tau''}
\mu_p(\sigma^{*,n-1}(\tau))\,d\tau\bigg]
\Bigg\}\,.
\end{multline*}
We observe that in this case
\begin{equation}
\label{mu}
\mu_p\big(\sigma^{*,n-1}(\tau)\big)\equiv
\mu_p\big(\sigma^{n-1}(0)\big),
\end{equation}
and hence
$$
\int\limits_{\tau'}^{\tau''}
\mu_p\big(\sigma^{*,n-1}(\tau)\big)\,d\tau=
\mu_p\big(\sigma^{n-1}(0)\big)(\tau''-\tau')\,.
$$
It follows that
\begin{multline}
\label{star}
I_2(-\tau',\tau')\le C_7 
\max\Bigg\{\Bigg[\int\limits_{-\tau''-1}^{-\tau''}d\tau\,
\bigg(\int\limits_{\sigma^{*,n-1}(\tau)}
|f(x)-C|^p\,d{\mathcal H}^{n-1}\bigg)^{1/(1-p)}\Bigg]^{1-p},
\\
\Bigg[\int\limits_{\tau''}^{\tau''+1}d\tau\,
\bigg(\int\limits_{\sigma^{*,n-1}(\tau)}
|f(x)-C|^p\,d{\mathcal H}^{n-1}\bigg)^{1/(1-p)}\Bigg]^{1-p}
\Bigg\}\\
\times
\exp\Big[-{\frac{\nu_1}{\nu_2}}\,\mu_p(\sigma^{n-1}(0))(\tau''-\tau')
\Big]\,.
\end{multline}
By letting $\beta,\tau''\to +\infty$ we obtain the following statement:

\begin{thm}
Fix a domain $D_0$ in ${\mathbb R}^{n-1}$ with compact and smooth boundary.
Let 
\[
D=D_0 \times \R = \{x \in {\mathbb R}^n : (x_1,\ldots,x_{n-1}) \in D_0\},
\]
and let $f\colon D \to \R$ be a generalized solution of
{\rm (\ref{eq2.19})} with {\rm (\ref{areq5})} and
{\rm (\ref{eq2.23})} satisfying
the boundary condition {\rm (\ref{areq3})} on $\partial D$.
If for a constant $C$ the right side of {\rm (\ref{star})} goes to $0$ as
$\tau''\to\infty$, then $f\equiv{\rm const}$ on $D$.
\end{thm}

Similarly for a solution $f$ of (\ref{eq2.19}) with (\ref{areq5}) and
(\ref{eq2.23}), satisfying the boundary condition (\ref{areeq3}) we may
write
\begin{multline}
\label{star2}
I_2(-\tau',\tau')\le C_7
\max\Bigg\{\Bigg[\int\limits_{-\tau''-1}^{-\tau''}d\tau\,
\bigg(\int\limits_{\sigma^{*,n-1}(\tau)}
|f(x)|^p\,d{\mathcal H}^{n-1}\bigg)^{1/(1-p)}\Bigg]^{1-p}, \\
\Bigg[\int\limits_{\tau''}^{\tau''+1}d\tau\,
\bigg(\int\limits_{\sigma^{*,n-a}(\tau)}
|f(x)|^p\,d{\mathcal H}^{n-1}\bigg)^{1/(1-p)}\Bigg]^{1-p}
\Bigg\}\,\\
\times
\max\Bigg\{\exp\bigg[-{\frac{\nu_1}{\nu_2}}\,\int\limits_{-\tau''}^{-\tau'}
\lambda_{p,Z^*_f(\tau)}^{1/p}(\sigma^{*,n-1}(\tau))\,d\tau\bigg],
\exp\bigg[-{\frac{\nu_1}{\nu_2}}\,\int\limits_{\tau'}^{\tau''}
\lambda_{p,Z^*_f(\tau)}^{1/p}(\sigma^{*,n-1}(\tau))\,d\tau\bigg] \Bigg\}\,.
\end{multline}
However here we do not have any identity similar to (\ref{mu}). We have:

\begin{thm}
Fix a domain $D_0$ in ${\mathbb R}^{n-1}$ with compact and smooth boundary.
Let 
\[
D=D_0 \times \R = \{x \in {\mathbb R}^n : (x_1,\ldots,x_{n-1}) \in D_0\},
\]
and let $f\colon D\to \R$ be a generalized
solution of {\rm (\ref{eq2.19})} with
{\rm (\ref{areq5})} and {\rm (\ref{eq2.23})} satisfying the boundary condition {\rm (\ref{areeq3})} on $\partial D$. If the right side of {\rm (\ref{star2})} tends to $0$ as $\tau''\to\infty$, then $f\equiv{\rm 0}$ on $\partial D$.
\end{thm}

If $f(x)=0$ everywhere on $\partial D$, then an identity similar to (\ref{mu}) holds in the 
following form:
\begin{equation}
\label{lambda}
\lambda_{p}^{1/p}(\sigma^{*,n-1}(\tau))\equiv
\lambda_{p}^{1/p}(\sigma^{n-1}(0))
\,.
\end{equation}
As above, we find
\begin{multline}
\label{star3}
I_2(-\tau',\tau')\le C_7
\max\Bigg\{\Bigg[\int\limits_{-\tau''-1}^{-\tau''}d\tau\,
\bigg(\int\limits_{\sigma^{*,n-1}(\tau)}
|f(x)|^p\,d{\mathcal H}^{n-1}\bigg)^{1/(1-p)}\Bigg]^{1-p},\\
\Bigg[\int\limits_{\tau''}^{\tau''+1}d\tau\,\bigg(\int\limits_{\sigma^{*,n-1}(\tau)}
|f(x)|^p\,d{\mathcal H}^{n-1}\bigg)^{1/(1-p)}\Bigg]^{1-p}
\Bigg\}\\
\times
\exp\bigg[-{\frac{\nu_1}{\nu_2}}\,\lambda_{p}^{1/p}(\sigma^{n-1}(0))(\tau''-\tau')
\bigg].
\end{multline}
Thus we obtain:

\begin{cor}
Fix a domain $D_0$ in ${\mathbb R}^{n-1}$ with compact and smooth boundary.
Let 
\[
D=D_0 \times \R = \{x \in {\mathbb R}^n : (x_1,\ldots,x_{n-1}) \in D_0\},
\]
and let $f\colon D\to \R$ be a generalized solution of
{\rm (\ref{eq2.19})} with {\rm (\ref{areq5})} and
{\rm (\ref{eq2.23})} satisfying
the boundary condition $f=0$ on $\partial D$.
If the right side of {\rm (\ref{star3})} tends to $0$ as $\tau''\to\infty$,
then $f\equiv{\rm const}$ on $\partial D$.
\end{cor}



\subsection*{Phragm\'en-Lindel\"of type theorems II}{}

We prove Phragm\'en-Lindel\"of type theorems for canonical domains of an
arbitrary form. Let $1\le k<n-1$. We consider a
domain
\[
D=D_0 \times \R^{n-k} = \{x \in {\mathbb R}^n : (x_1,\ldots,x_{k}) \in D_0\},
\]
where $D_0$ is a domain in ${\mathbb R}^{k}$ with compact and smooth boundary.
Let $f$ be a generalized solution of (\ref{eq2.19}) with (\ref{areq5})
and (\ref{eq2.23}) satisfying
the boundary condition (\ref{areq3}) on $\partial D$.

Fix $\tau_0>0$. Let $\tau_0<\tau'<\tau''<\infty\,.$
By (\ref{eq203}) we may write
$$
\int\limits_{D^k_{0,\tau'}}
|\nabla f|^p\,d{\mathcal H}^n\\
\le C_8
\Bigg[\int\limits_{\tau'}^{\tau''}d\tau\,\bigg(\int\limits_{\sigma^k(\tau)}
|f(x)-C|^p\,d{\mathcal H}^{n-1}\bigg)^{1/(1-p)}\Bigg]^{1-p}
\,,
$$
where $C_8=C_7/2$.
As in (\ref{eqSVP2}) we obtain from (\ref{eeqIlI}) the estimate
$$
\int\limits_{D^k_{0,\tau_0}}
|\nabla f|^p\,d{\mathcal H}^n
\le \int\limits_{D^k_{0,\tau'}}
|\nabla f|^p\,d{\mathcal H}^n
\,
\exp\bigg[-{\frac{\nu_1}{\nu_2}}\,\int\limits_{\tau_0}^{\tau'}
\mu_p(\sigma^k(\tau))\,d\tau\bigg]\,.
$$
By combining these inequalities we obtain
\begin{multline}
\label{eq7.12}
\int\limits_{D^k_{0,\tau_0}}
|\nabla f|^p\,d{\mathcal H}^n
\le C_8
\Bigg[\int\limits_{\tau'}^{\tau''}d\tau\,\bigg(\int\limits_{\sigma^k(\tau)}
|f(x)-C|^p\,d{\mathcal H}^{n-1}\bigg)^{1/(1-p)}\Bigg]^{1-p}\\
\times
\exp\bigg[-{\frac{\nu_1}{\nu_2}}\,\int\limits_{\tau_0}^{\tau'}
\mu_p(\sigma^k(\tau))\,d\tau\bigg]\,.
\end{multline}
The inequality (\ref{eq7.12}) holds for arbitrary
constant $C$ and every $\tau''>\tau'$. Thus the following statement
holds:

\begin{thm}
Let $f\colon D \to \R$ be a generalized solution
of {\rm (\ref{eq2.19})} with
{\rm (\ref{areq5})} and
{\rm (\ref{eq2.23})} satisfying the boundary condition {\rm (\ref{areq3})} on
$\partial D$, $1\le k<n-1$. If for a constant $C$ the right side of
{\rm (\ref{eq7.12})} tends to $0$ as $\tau',\,\tau''\to +\infty$, then
$f\equiv {\rm const}$ on $D$.
\end{thm}

If $f$ satisfies (\ref{eq2.19}) with (\ref{areq5}),
(\ref{eq2.23}) and the boundary condition (\ref{areeq3}) on
$\partial D$,
then we have
\begin{multline}
\label{eq10.39}
\int\limits_{D^k_{0,\tau_0}}
|\nabla f|^p\,d{\mathcal H}^n
\le
C_8
\Bigg[\int\limits_{\tau'}^{\tau''}d\tau\,\bigg(\int\limits_{\sigma^k(\tau)}
|f(x)|^p\,d{\mathcal H}^{n-1}\bigg)^{1/(1-p)}\Bigg]^{1-p}\\
\times
\exp\bigg[-{\frac{\nu_1}{\nu_2}}\,\int\limits_{\tau_0}^{\tau'}
\lambda_{p,Z_f(\tau)}^{1/p}(\sigma^k(\tau))\,d\tau\bigg]\,.
\end{multline}
We obtain:

\begin{thm}
Fix a domain $D_0$ in ${\mathbb R}^{k}$, where $1\le k<n-1$, with compact and smooth boundary. Let 
\[
D=D_0 \times \R^{n-k} = \{x \in {\mathbb R}^n : (x_1,\ldots,x_{k}) \in D_0\},
\]
and let $f\colon D\to \R$ be a generalized solution of {\rm (\ref{eq2.19})}  with
{\rm (\ref{areq5})} and {\rm (\ref{eq2.23})} satisfying the boundary conditions {\rm (\ref{areeq3})} on $\partial D$. If for a constant $C$ the right side of {\rm (\ref{eq10.39})} tends to $0$ as $\tau',\,\tau''\to +\infty$, then $f\equiv 0$ on $D$.
\end{thm}


\end{document}